I. Sh. Jabbarov

On Ergodic Hypothesis.

## 1. Introduction

In the theory of stochastic processes The Ergodic Hypothesis of Boltzman (see [6, 7, 8, 14]) is well known. This hypothesis asserts that the average time of a finding of physical systems in some domain of phase space is equal to the relative measure of given domain (see [6], p. 701, [7], p. 522). However, the movement trajectory, in this case, depends on an initial position of a point, and the average on time undertaken irrespectively of it. The exact mathematical formulation of this hypothesis is expressed by a limit relation of a type

$$\lim_{T \to \infty} \frac{1}{T} \int_0^T f(T_t(x)) dt = \int_\Omega f(\theta) \mu(d\theta), \quad (1)$$

where $T_t$ some semigroup of mappings of phase space (see [6], [7]), $x$ is, generally speaking, any point of the phase space. We will, at the expense of generality, to consider the problem in such statement. In [2, p. 457] this case is considered for continuous functions $f$ in "a narrowed ergodic" case.

Let's notice that the Ergodic hypothesis, generally, is not true for the function (the example in [8] see, p. 536), accepting value 1 on a given trajectory, and 0, on other points of the phase space. Analogically, it is possible to construct a similar example in the infinite dimensional case. This is clear from the definition of product Lebesgue measure.

In the works [6, 7, 8, 14] some have been proved ergodic theorems in which the relation (1) is proved for almost all points $x$ of the phase space.

In the present work we show that, if as a phase space to consider the cube $\Omega = [0,1] \times [0,1] \times \cdots$, it is possible to enter such a dynamical system for which the hypothesis mentioned above holds for a certain class of functions and any point $x$.

Let's consider the dynamical system defined by a semigroup of mappings: $x \to \varphi_t(x)$, where $\varphi_t(x) = \{x + t\overline{\Lambda}\} = (\{x_n + t\lambda_n\})$, $x \in \Omega$, (the symbol $\{\cdot\}$ means a fractional part), $\overline{\Lambda} = (\lambda_n)$, $\lambda_n \to \infty$ is a sequence of positive numbers. Let $f$ be a function defined in $\Omega$, $f \in L_2(\Omega, \mu_0)$ where $\mu_0$ designates the measures in the $\Omega$ introduced in the works [9, 11, 12, 13]. Function $f = f(\overline{\theta}) = f(\theta_1, \theta_2, ...)$ can be expanded into the Fourier series

$$f \sim \sum a(\overline{m}) e^{2\pi i (\overline{m}, \overline{\theta})},$$

where $\overline{\theta} \in \Omega$, $\overline{m} = (m_1, m_2, ...)$ and $m_k = 0$ for all $m_k$, with exception of finite number of them. Let for any natural $r \in N$ the multiple series

$$\sum_{(m_1, m_2, ..., m_r) \in Z^r} |a(m_1, m_2, ..., m_r, 0, ...)|$$

converges (i.e. the corresponding subseries of given Fourier series converges absolutely). We will say that the function $f$ belongs to the class $J_\infty$ if this condition satisfied for any natural $r$. Below we use the notions and designations defined in the works [9, 11, 12, 13].

**Definition 1.** *Let $\sigma : N \to N$ be any one to one mapping of the set of natural numbers. If for any $n > m$ there is a natural number $m$ such that $\sigma(n) = n$, then we call $\sigma$ a finite permutation. A subset $A \subset \Omega$ is called to be finite-symmetrical if for any element $\theta = (\theta_n) \in A$ and any finite permutation $\sigma$ one has $\sigma\theta = (\theta_{\sigma(n)}) \in A$.*

Let $\overline{\omega} \in \Omega$ and $\Sigma(\overline{\omega}) = \{\sigma\overline{\omega} \mid \sigma \in \Sigma\}$. It is a countable set. We will designate $\Sigma'(\overline{\omega})$ the set of all limit points of the sequence $\Sigma(\overline{\omega})$. For every real $t$ we write $\{t\overline{\Lambda}\} = (\{t\lambda_n\})$. Let

$$f_r(\overline{\theta}) = \sum a(\overline{m}) e^{2\pi i (\overline{m}, \overline{\theta})} = \sum_{(m_1, ..., m_r) \in Z^r} a(m_1, ..., m_r, 0, ...) e^{2\pi i (\overline{m}, \overline{\theta})}.$$

The following theorem holds:

**Theorem 1.** *There is a sequence of natural numbers $(r_k)$ and a subset $\Omega_0 \subset \Omega$ of full measure such that the series $\sum_{k \geq 2} |f_{r_k}(\overline{\theta}) - f_{r_{k-1}}(\overline{\theta})|$ converges everywhere in $\Omega_0$.*

Thus, the sequence $(f_{r_k})$ converges in $\Omega_0$ to some function $\overline{f}(\theta)$. The relation $\overline{f}(\theta) = f(\theta)$ is fair almoct everywhere in $\Omega$ and we will prove the relation (1) at first for the function $\overline{f}(\theta) = f(\theta)$.

**Theorem 2.** *For any function $f \in J_\infty$ the following relation holds*

$$\lim_{T \to \infty} \frac{1}{T} \int_0^T \overline{f}(\overline{u} + (\{t\lambda_n\})) dt = \int_\Omega \overline{f}(\theta) \mu_0(d\theta) = \int_\Omega f(\theta) \mu_0(d\theta),$$

*for any point $\overline{u} \in \Omega$.*

The increasing sequence of linearly independent frequencies $\lambda_n$ plays here an essential role. In H. Bohr's works [3, 4] averages of a kind (1) on the basis of Croneker's theorem on the uniform distribution (mod1) of some curves in the multidimensional unite cube are studied. The theorem of Croneker (see [18, p. 301]) states:

*Let $\alpha_1, \alpha_2, ..., \alpha_N$ be real numbers linearly independent over the field of rational numbers, $\gamma$ be a subdomain of $N$ – dimensional unite cube with the volume $\Gamma$ in Jordan meaning. Let further, $I_\gamma(T)$ is a measure of a set of such $t \in (0,T)$ for which $(\alpha_1 t, \alpha_2 t, ..., \alpha_N t) \in \gamma \pmod 1$. Then*

$$\lim_{T \to \infty} \frac{I_\gamma(T)}{T} = \Gamma.$$

Later the theory of uniform distribution of curves has been generalized and deeply studied by many other authors (see [15]). Generalization of this theory to the curves given by functions integrable in the Lebesgue sence meets great difficulties, if as a subdomain $\gamma$ to take a subset measurable in the Lebesgue sence (see [15 p. 96]).

The question on continuous uniform distribution of (mod1) curve $\bar{f}(\{t\Lambda\})$ is connected with an existence of limit of a kind

$$\lim_{T \to \infty} \frac{1}{T} \int_0^T w(\{\bar{f}(t)\}) dt$$

where $w(x)$ is continuous in the cube $[0,1]^n$ (see [9]). The similar limit is studied in the theory of dynamical systems where the question is connected with the metric transitivity. However, existence of the specified limit is proved only for almost all shifts of trajectory in the phase space (see [6]).

## 2. Basic auxiliary lemmas.

**Lemma 1**. *Let $A \subset \Omega$ be a finite-symmetrical subset of zero measure and $\Lambda = (\lambda_n)$ is an unbounded, monotonically increasing sequence of positive real numbers any finite subfamily of elements of which is linearly independent over the field of rational numbers. Let $B \supset A$ be any open, in the Tychonoff metric, subset with $\mu_0(B) < \varepsilon$,*

$$E_0 = \{0 \le t \le 1 \mid \{t\Lambda\} \in A \wedge \Sigma'\{t\overline{\Lambda}\} \subset B\}.$$

*Then, $m(E_0) \le c_0 \varepsilon$ where $c_0 > 0$ is an absolute constant, $m$ designates the Lebesgue measure.*

Proof of this lemma is given in [9, 11].

**Lemma 2**. *Let a curve $\gamma(t)$ to be uniformly distributed (mod 1) in the space $\mathbf{R}^n$. Let D be a closed subdomain of a unite cube measurable in Jordan meaning, and $\Phi$ be a family of continuous* *functions defined in D. If $\Phi$ is a uniformly bounded and equicontinuouc family then uniformly by $f \in \Phi$ the relation below holds*

$$\lim_{T \to \infty} T^{-1} \int\limits_{\{\gamma(t)\} \in D} f(\{\gamma(t)\}) dt = \int\limits_D f dx_1 ... dx_N,$$

*where on the left part of the equality the integration is taken over such $t \in (0,T)$ for which $\gamma(t) \in D(\mathrm{mod}\,1)$ and $\{\gamma(t)\} = (\{\gamma_1(t)\},...,\{\gamma_N(t)\})$.*

In particular, the lemma 2 is fair for any continuous function in $D$. Proof of the lemma 2 can be found in [19, p. 348].

**Lemma 3.** *Under the conditions of the lemma 1, for any continuous in the closed ball $\overline{B} = \overline{B}(\theta_0, r) \subset \Omega$ function $f$, the relatin*

$$\lim_{T \to \infty} T^{-1} \int\limits_{0, \{t\overline{\Lambda}\} \in B}^{T} f(\{t\overline{\Lambda}\}) dt = \int\limits_B f(\overline{\theta}) \mu_0(d\overline{\theta})$$

*holds*.

*Proof.* The lemma 3 follows from the lemma 2, if to apply it, at first, to the suitable projections (see[9, 13]) $B_N(\overline{\theta}_0, r - \varepsilon)$ and $B_N(\overline{\theta}_0, r)$, and then to pass to the limit at $\varepsilon \to 0$, since

$$\int\limits_{\{t\overline{\Lambda}\} \in B_N(\overline{\theta}_0, r-\varepsilon)} f(\{t\overline{\Lambda}\}) dt \leq \int\limits_{\{t\overline{\Lambda}\} \in B} f(\{t\overline{\Lambda}\}) dt \leq \int\limits_{\{t\overline{\Lambda}\} \in B_N(\overline{\theta}_0, r)} f(\{t\overline{\Lambda}\}) dt.$$

The lemma 3 is proved.

**Lemma 4.** *Let $A \subset \Omega$ be a closed measurable subset, $f$ - continuous in A function in the Tychonoff metrics, (i.e. If $\overline{\theta} \in A$ and $\overline{\theta}_m \to \overline{\theta}, \overline{\theta}_m \in \Omega$ then $f(\overline{\theta}_m) \to f(\overline{\theta})$ as $m \to \infty$). Then the relation below holds*

$$\lim_{T \to \infty} T^{-1} \int\limits_{0, \{t\overline{\Lambda}\} \in A}^{T} f(\{t\overline{\Lambda}\}) dt = \int\limits_A f(\overline{\theta}) \mu_0(d\overline{\theta}).$$

*Proof.* Since $A$ is closed, there is a finite number of open balls the union $E$ of which contains $A$ (see [1, p. 221]). Thus, for any $\varepsilon > 0$, it is possible to assume that $\mu_0(A) \leq \mu_0(E) \leq \mu_0(A) + \varepsilon$. Let $|f|$ has a maximum of the modulus $M$ in $A$. Every ball $B$ we can include into some set of view $S \times \Omega$, where $S$ is a projection of the ball $B$ into the subspace of first $N$ coordinate axes, and $N$ can be chosen such that $\mu_0(B) \leq \mu_0(S) + \varepsilon$ (see[11]). Since, the $A$ is overlapped by the finite union of open balls then the union of these balls also can be included into the set of view $S \times \Omega$ with the same property for some natural $N$. So, we have $\mu_0(A) \leq \mu_0(S) + 2\varepsilon$. Under the theorem on continuous continuation (see [16, p. 323]) the function $f_S(\bar{x}) = \int_\Omega f(\bar{x},\bar{\theta})\mu_0(d\theta)$ $(\bar{x},\bar{\theta}) \in A$ has an extension to $R^N$ with the maximum of the modulus $M$. Owing to a continuity of the extended function on the closed set $\bar{S} \times \Omega$, this function will be uniformly continuous (in the Tychonoff metric). Therefore, there exist $\delta > 0$ such that $|f(\bar{\theta}) - f(\bar{\theta}')| < \varepsilon$ if $|\bar{\theta} - \bar{\theta}'| < \delta$. Since

$$\sum_{n>N} e^{1-n}|\bar{\theta}_n - \bar{\theta}'_n| \leq e^{1-N},$$

then it could be found $N_0$ such that for any $N > N_0$ and $\bar{x} \in S \subset R^N$ the relation

$$\left| f(\bar{x},\bar{\theta}_0) - \int_\Omega f(\bar{x},\bar{\theta})\mu_0(d\bar{\theta}) \right| \leq \varepsilon$$

is satisfied at any point $\bar{\theta}_0 \in \Omega$. Applying the lemma 2, we find:

$$\lim_{T\to\infty} T^{-1} \int_{0,\{t\bar{\Lambda}\}\in A}^T f(\{t\bar{\Lambda}\})dt \leq \int_{\bar{S}}\int_\Omega f(\bar{x},\bar{\theta})d\bar{x}\mu_0(d\bar{\theta}) = \int_A f(\bar{\theta})\mu_0(d\bar{\theta}) + O(\varepsilon).$$

Since the set $A$ is measurable then the set $A' = (\Omega \setminus \overline{A}) \cup \partial A$ is a closed measurable set. So, the reasoning above give the relation

$$\lim_{T\to\infty} T^{-1} \int_{0,\{t\bar{\Lambda}\}\in A'}^T f(\{t\bar{\Lambda}\})dt \leq \int_{A'} f(\bar{\theta})\mu_0(d\bar{\theta}) + O(\varepsilon)$$

Therefore,

$$\int_{A'} f(\bar{\theta})\mu_0(d\bar{\theta}) - O(\varepsilon) \leq \lim_{T\to\infty} T^{-1} \int_{0,\{t\bar{\Lambda}\}\in A}^T f(\{t\bar{\Lambda}\})dt \leq \int_A f(\bar{\theta})\mu_0(d\bar{\theta}) + O(\varepsilon)$$

Then tending $\varepsilon$ to zero, we reseive:

$$\lim_{T\to\infty} T^{-1} \int_{0,\{t\bar{\Lambda}\}\in A}^T f(\{t\bar{\Lambda}\})dt = \int_A f(\bar{\theta})\mu_0(d\bar{\theta}).$$

The lemma 4 is proved.

### 3. Proofs of the theorems.

**Proof of the theorem 1.** We have:

$$\int_\Omega \sum_{k\geq 2} |f_{r_k}(\bar{\theta}) - f_{r_{k-1}}(\bar{\theta})|\mu_0(d\bar{\theta}) \leq \sum_{k\geq 2}\left( \int_\Omega |f_{r_k}(\bar{\theta}) - f_{r_{k-1}}(\bar{\theta})|^2 \mu_0(d\bar{\theta}) \right)^{1/2}. \qquad (2)$$

As $f \in L_2(\Omega,\mu_0)$, then the series

$$\sum_{\bar{m}\in Z^\infty} |a(\bar{m})|^2$$

converges. So, it could be found a sequence $(r_k)$ of natural numbers such that

$$\int_\Omega |f_{r_k}(\bar\theta) - f_{r_{k-1}}(\bar\theta)|^2 \mu_0(d\bar\theta) \leq 2^{-2k}. \qquad (3)$$

Consequently, the series on the right side of the inequality (2) is convergent. Then, there is a subset $\Omega_0 \subset \Omega$ of a full measure such that the series $\sum_{k\geq 2}|f_{r_k}(\bar\theta) - f_{r_{k-1}}(\bar\theta)|$ converges pointwisely in $\Omega_0$. Then, in the $\Omega_0$ such a function $\bar f$ is defined that $f_{r_k} \to \bar f$. The theorem 1 is proved.

**Proof of the theorem 2.** By the definition of the functions $f_{r_k}$, and condition od absolute convergense of Fourier series for them, the divergence of a series defining $\bar f$ does not change after acting of any finite permutation to the point $\bar\theta$, i.e. the set $\Omega \setminus \Omega_0$ where the series $\sum_{k\geq 2}|f_{r_k}(\bar\theta) - f_{r_{k-1}}(\bar\theta)|$ diverges is finite-symmetrical. Under the Egoroff's theorem, this series converges almost uniformly (see [1, p.111]) in $\Omega$. We need in more exact characteristics of a set where the considered seties defines a continuous function. For receive it recall the construc-tion of a proof of Egoroff's theorem (see [13]). Let $(\varepsilon_k)$ be a sequence of positive numbers, $\varepsilon_k \to 0$. Consider in the set $\Omega_0$ of convergence of the series $\sum_{m\geq 2}|f_{r_m}(\bar\theta) - f_{r_{m-1}}(\bar\theta)|$ the sequence $g_n = \sum_{m>n}|f_{r_m} - f_{r_{m-1}}|$. Denote $S_{k,l} = \{\bar\theta \in \Omega_0 \mid n \geq k \wedge g_n < \varepsilon_l\}$, i. e. $S_{k,l}$ is the set of such $\bar\theta \in \Omega_0$ for which the inequality $g_n(\bar\theta) < \varepsilon_l$ is satisfied for every $n \geq k$. We have $S_{1,l} \subset S_{2,l} \subset \cdots$ and $\Omega_0 = \bigcup_{n\geq 1} S_{n,l}$ for each natural number $l \geq 1$. Let $\delta > 0$ be any positive number. One can define a natural number $n(l)$ for which $\mu_0(\Omega_0 \setminus S_{n(l),l}) < 2^{-l}\delta$. Denote now $U_k = \bigcap_{l=k}^\infty S_{n(l),l}$. In every set $U_k$ the series $\sum_m |f_{r_m} - f_{r_{m-1}}|$ converges uniformly and

$$\mu_0(\Omega_0 \setminus U_k) \leq \sum_{l\geq k}\mu_0(\Omega_0 \setminus S_{n(l),l}) < 2^{1-k}\delta.$$

Further, $U_1 \subset U_2 \subset \cdots$. Let $U = \bigcup_{k\geq 1}U_k$. So, $\mu_0(\Omega_0 \setminus U) = 0$ and if $\bar\theta \in \Omega_0 \setminus U$ then $\bar\theta \notin U_k$ for each $k$. Hence, for each natural $k$ one can find $l_k$ so that $\bar\theta \notin S_{n(l_k),l_k}$. We can suppose that $n(l_k) \to \infty$ as $k \to \infty$ (if else, then $n(l_k) \leq L$ for some natural $L$, and therefore, $g_n < 2^{1-k}\delta$ for some $n \leq L$ and all $k$; then, $g_n = 0$ from which one deduces that the considered Fourier series is absolutely convergent and the statement of the theorem is true). Now the relation $\bar{\theta} \notin S_{n(l_k),l_k}$ shows that there exist a sequence of natural numbers $(j_k)$ such that $j_k \geq n(l_k)$ and $g_{j_k}(\bar{\theta}) \geq \varepsilon_{l_k}$. If $\sigma$ is a given finite permutation then on the all great enough values of $k$ the components of $\bar{\theta}$ invariable after of action of $\sigma$, so $g_{j_k}(\sigma\bar{\theta}) \geq \varepsilon_{l_k}$. It means $\sigma\bar{\theta} \notin S_{n(l_k),l_k}$ for all great enough $k$, so that $\sigma\theta \notin U_k$. Therefore, $\sigma\bar{\theta} \in \Omega_0 \setminus U$ and this relation implies that the set $\Omega_1 = \Omega_0 \setminus U$ is finite symmetrical. *Moreover, if we should exchange the first $n(l_k)-1$ components of the point $\bar{\theta} \in \Omega_1$ by any numbers from the interval* [0,1]*, we get again some point $\bar{\theta}' \in \Omega_1$.* Since, $l_k \to \infty$ then the told above is fair for any first $n$ natural components. Denoting $\Omega' = \Omega_1 \cup (\Omega \setminus \Omega_0)$, we can state that for each $\bar{\theta} \in \Omega', \bar{\theta} = (\theta_1,\theta_2,...)$ we have $(\alpha_1,...,\alpha_n,\theta_{n+1},\theta_{n+2},...) \in \Omega'$, also, when $(\alpha_1,...,\alpha_n) \in [0,1]^n$ is any point and $n$ is a natural number. The needed characteristics are gotten.

Consider now the curve $(\{t\bar{\Lambda}\})$ in some open covering $B$ containing $\Omega'$. In agree with the theorem of Egoroff $B$ could be set such that $\mu_0(B) \leq \varepsilon$, for given $\varepsilon > 0$, and the series $\sum_{k \geq 2} |f_{r_k}(\bar{\theta}) - f_{r_{k-1}}(\bar{\theta})|$ converges uniformly in the $\Omega \setminus B$. The set $B$ can be represented as a union of open balls: $B = \bigcup_r B_r$. For every natural $n$ we define the set $\Sigma'_n(\bar{\omega}), \bar{\omega} \in \Omega$ ([11]) as a set of all limit points of the sequence

$$\Sigma_n(\omega) = \{\sigma\omega \mid \sigma \in \Sigma \wedge \sigma(1) = 1, \wedge \cdots \wedge \sigma(n) = n\}.$$

Let

$$D^{(n)} = \{t \mid \{t\bar{\Lambda}\} \in \Omega_1 \wedge \Sigma'_n(\{t\bar{\Lambda}\}) \subset \bigcup_r B_r\}$$

we have $D^{(n)} \subset D^{(n+1)}$. Denote by $D$ the outer limit set $D = \bigcup D^{(n)}$. Since the conditions of the lemma 1 are satisfied then we deduce $m(D^{(n)}) \leq c_0\varepsilon$. Thus, we have $m(D) \leq c_0\varepsilon$. Hence, if $t \notin D$, taking $n = r_k, k = 1,2,...$, we can find such limit point $\bar{\omega}_k \in \Omega \setminus \bigcup_{r=1}^{\infty} B_r$ of the sequence $\Sigma_n(\{t\bar{\Lambda}\})$, for which the series

$$\sum_{m \geq 2} |f_{r_m}(\omega_k) - f_{r_{m-1}}(\omega_k)|$$

converges for all values of $k$. As the set $\Omega \setminus \bigcup_{r=1}^{\infty} B_r$ is closed, the limit point $(\{t\bar{\Lambda}\})$ of the sequence $(\bar{\omega}_k)$ will belong to the set $\Omega \setminus \bigcup_{r=1}^{\infty} B_r$. Therefore, the series

$$\sum_{m \geq 2} \left| f_{r_m}(\{t\overline{\Lambda}\}) - f_{r_{m-1}}(\{t\overline{\Lambda}\}) \right|$$

converges for all $t \in E \setminus D$, i.e. $E \subset D$ where $E$ is a subset of such real numbers $t, 0 \leq t \leq 1$ that the series $\sum_{m \geq 2} \left| f_{r_m}(\{t\overline{\Lambda}\}) - f_{r_{m-1}}(\{t\overline{\Lambda}\}) \right|$ is divergent. Then, for all $t$, with exception of values from some set of a measure, not exceeding $2c_0\varepsilon$, the last series converges. Owing to randomness of $\varepsilon$, this result shows a convergence of the series of the lemma 2 for almost all considered $t$. Moreover, $\mu_0(D) \leq \varepsilon$. It is clear that this is true olso for almost all real $t$.

Let $g(\overline{\theta}) = |f_{r_1}(\overline{\theta})| + \sum_{m \geq 2} |f_{r_m}(\overline{\theta}) - f_{r_{m-1}}(\overline{\theta})|$. Consider now the sequence $(g(\sigma\overline{\theta}))_{\sigma \in \Sigma}$ and define the function $\inf\{g(\sigma\overline{\theta}) \mid \sigma \in \Sigma\} = \rho(\overline{\theta})$ for each considered $\overline{\theta}$. This function is a measurable function. Define now subsets

$$\Omega_k = \{\overline{\theta} \in B \mid \rho(\overline{\theta}) \geq 2^k\}; \quad E_k = \{t, 0 \leq t \leq 1 \mid \{t\overline{\Lambda}\} \in \Omega_k\};$$

$$E_0 = \{t \mid 0 \leq t \leq 1 \wedge \{t\overline{\Lambda}\} \in D\}.$$

Since for each $\overline{\theta} \in \Omega_k$ we have $g(\overline{\theta}) \geq 2^k$, then from (2-3) we deduce

$$2^{2k} \mu_0(\Omega_k) < \int_\Omega g^2(\overline{\theta}) \mu_0(d\overline{\theta}) \leq \sum_{m=1}^\infty \frac{1}{m^2} \sum_{m=2}^\infty m^2 \int_\Omega |f_{r_m} - f_{r_{m-1}}|^2 \mu_0(d\overline{\theta}) + C, \quad (4)$$

where $C$ is a constant. So, we have $\mu_0(\Omega_k) \leq b \cdot 2^{-2k}$ for some positive constant $b$. We can estimate the Lebesgue measure of a set (which is placed in the closed and finite-symmetrical subset) $E_k$, by using of any finite permutation $\sigma$ satisfying the conditions of the lemma 4 of [13] as a value

$$m(E_k) \leq c_0 b 2^{-2k} \quad (5)$$

with some suitable constant $c_0 > 0$.

Consider now the $k$ such that $2^{-k} > \varepsilon$. To every point $\overline{\theta} \in B$, on which $\rho(\overline{\theta}) \geq 2^k$, we can put in correspondence some finite permutation $\sigma$ satisfying the condition

$$g(\sigma\overline{\theta}) < 2^{k+1}$$

in some neighborhood $W(\overline{\theta})$ of the point $\overline{\theta}$. If we define the compact subset $M = \{\overline{\theta} \in \Omega \mid g(\overline{\theta}) \leq 2\varepsilon^{-1}\}$ then $M$ can be overlapped by the union of all open neighborhoods $W(\overline{\theta})$ of the points $\overline{\theta}$. Therefore, there is only finite family of neighborhoods the union of which contains the set $M$. Since, the function $g(\sigma\overline{\theta})$ is continuous in the union $(\Omega \setminus B) \cup M$, for each permutation $\sigma$ (this permutation is defined by taken neighborhood $W(\overline{\theta})$), then in the set $(\Omega \setminus B) \cup M$ the lemma

3 is applicable:

$$\lim_{T\to\infty}\frac{1}{T}\int_{\{t\overline{\Lambda}\}\in(\Omega\setminus B)\cup M} g(\sigma\{t\overline{\Lambda}\})dt = \int_{(\Omega\setminus B)\cup M} g(\sigma\overline{\theta})\mu_0(d\overline{\theta}). \quad (6)$$

Consider now the set $B\setminus M$ which is overlapped by the union $(\Omega\setminus\Omega_0)\cup\left(\bigcup\Omega_k\right)$ taken over $k\geq 1$ satisfying the constraints $2^{-k}\leq\varepsilon$. Let $T_0$ be taken so that, for every $T>T_0$ and $\eta>0$

$$\left|\frac{1}{T}\int_{\{t\overline{\Lambda}\}\in(\Omega\setminus B)\cup M}\overline{f}(\{t\overline{\Lambda}\})dt - \int_{(\Omega\setminus B)\cup M}\overline{f}(\overline{\theta})\mu_0(d\overline{\theta})\right|<\eta.$$

We have

$$\int_0^T \overline{f}(\{t\overline{\Lambda}\})dt = \int_{\{t\overline{\Lambda}\}\in\Omega_1}\overline{f}(\{t\overline{\Lambda}\})dt + \sum_{m=2}^{\infty}\int_{\{t\overline{\Lambda}\}\in\Omega_m\setminus\Omega_{m-1}}\overline{f}(\{t\overline{\Lambda}\})dt$$

The integrals on the right side of the last equality are of the type of integtals taken over the subsets $D$ and $E_m\setminus E_{m-1}$, respectively (note that on sufficiently small values of $\varepsilon$ a subsets $E_m$ for small natural $m$ are empty). Then from (5) it follows that

$$\left|\sum_{n\leq T}\int_{\{t\overline{\Lambda}\}\in E_m\setminus E_{m-1},\,n\leq t\leq n+1} g(\{t\overline{\Lambda}\})dt\right|<4c_0 b\varepsilon(T+1). \quad (7)$$

Therefore, for every $T>T_0$

$$\left|\frac{1}{T}\int_{\{t\overline{\Lambda}\}\in B\setminus M,\,t\leq T} g(\{t\overline{\Lambda}\})dt\right|<5c_0 b\varepsilon.$$

Let, now, $T>0$ to be great enough, so that (5) is satisfied with an error not exceeding $\varepsilon$. We will consider the integral

$$T^{-1}\int_0^T \overline{f}(\{t\overline{\Lambda}\})dt = T^{-1}\int_{\{t\overline{\Lambda}\}\in(\Omega\setminus B)\cup M}\overline{f}(\{t\overline{\Lambda}\})dt + T^{-1}\int_{t\in B\setminus M}\overline{f}(\{t\overline{\Lambda}\})dt. \quad (8)$$

By the lemma 4

$$T^{-1}\int_{\{t\overline{\Lambda}\}\in(\Omega\setminus B)\cup M}\overline{f}(\{t\overline{\Lambda}\})dt = \int_{(\Omega\setminus B)\cup M} f(\overline{\theta})\mu_0(d\overline{\theta}) + o_T(1) = \int_\Omega f(\overline{\theta})\mu_0(d\overline{\theta}) + O(\sqrt{\varepsilon}), \quad (9)$$

since

$$\int_B g(\overline{\theta})\mu_0(d\overline{\theta})\leq\sqrt{\mu_0(B)}\left(\int_\Omega g^2(\overline{\theta})\mu_0(d\overline{\theta})\right)^{1/2}.$$

Applying the relations (7-9) and tending $\varepsilon\to 0$ we get a sutable result. The theorem 2 is proved for the case $\overline{u}=0$. General case can be considered by the same way.